\documentclass[11pt]{article}

\usepackage{amsthm}
\usepackage{amsmath}
\usepackage{amssymb}

\usepackage{graphicx}

\newtheorem*{theorem}{Theorem}
\newtheorem{ntheorem}{Theorem}
\newtheorem*{proposition}{Proposition}
\newtheorem*{corollary}{Corollary}
\newtheorem{nproposition}{Proposition}

\newcommand{\x}{{\boldsymbol{x}}}

\newcommand{\Chi}{\overline{\chi}}
\newcommand{\1}{\overline{1}}
\newcommand{\2}{\overline{2}}
\newcommand{\3}{\overline{3}}
\newcommand{\4}{\overline{4}}
\newcommand{\5}{\overline{5}}

\renewcommand{\j}{\overline{j}}

\newcommand{\n}{\overline{n}}
\newcommand{\Q}{\mathbb{Q}}
\newcommand{\C}{\mathbb{C}}
\newcommand{\N}{\mathbb{N}}
\newcommand{\R}{\mathbb{R}}
\newcommand{\A}{{\cal A}}
\renewcommand{\o}{{\it o}}

\newcommand{\+}{\oplus}

\renewcommand{\x}{\odot}

\renewcommand{\Im}{\mathrm{Im}}
\newcommand{\first}{\mathrm{first}}
\newcommand{\last}{\mathrm{last}}
\newcommand{\iin}{\mathrm{in}}
\renewcommand{\span}{\mathrm{span\,}}
\newcommand{\smallcells}{\mathrm{small\,}}
\newcommand{\largecells}{\mathrm{large\,}}
\newcommand{\trop}{\mathrm{trop\,}}
\newcommand{\troppos}{\mathrm{trop^+\,}}

\begin{document}
\title{A tropical morphism related to the hyperplane arrangement of the complete bipartite graph}
\author{Federico Ardila \footnote{Mathematical Sciences Research
Institute \, -- \, \textsf{federico@msri.org}}}
\date{}
\maketitle

\begin{abstract}
We undertake a combinatorial study of the piecewise linear map $g:
\R^{2(m+n)} \rightarrow \R^{m \times n}$ which assigns to the four
vectors $a, A$ in $\R^m$ and $b, B$ in $\R^n$ the $m \times n$
matrix given by $g_{ij}=\min(a_i+b_j, A_i+B_j)$. This map arises
naturally in Pachter and Sturmfels's work on the tropical geometry
of statistical models. The image of $g$ has been a subject of
recent interest; it is the positive part of the tropical algebraic
variety which parameterizes $n$-tuples of points on a tropical
line in $m$-space.

The domains of linearity of $g$ are the regions of the real
hyperplane arrangement $\A_{m,n}$, corresponding to the complete
bipartite graph $K_{m,n}$. We explain how the images of (some of)
the regions provide two polyhedral subdivisions of the image of
$g$, one of which is a refinement of the other. The finer
subdivision is particularly nice enumeratively: it has $2{m
\choose 2}{n \choose 2}r_{m-2,n-2}$ maximum-dimensional cells,
where $r_{m-2,n-2}$ is the number of regions of the arrangement
$\A_{m-2,n-2}$.

\end{abstract}

\section{Introduction}\label{sec:intro}

0. The goal of this paper is to undertake a combinatorial study of
the piecewise linear map $g: \R^{2(m+n)} \rightarrow \R^{m \times
n}$ given by
\[
g(a,A,b,B)_{ij} = \min(a_i + b_j, A_i + B_j) \qquad \mathrm{for}
\,\,\,\, 1 \leq i \leq m, 1 \leq j \leq n,
\]
where $a$ and $A$ denote vectors in $\R^m$ and $b$ and $B$ denote
vectors in $\R^n$.

The paper is organized as follows. In Section \ref{sec:faces}, we
describe and enumerate the domains of linearity of the map $g$;
they are precisely the faces of the arrangement $\A_{m,n}$,
corresponding to the complete bipartite graph $K_{m,n}$. In
Section \ref{sec:images} we describe the images $g(R)$, as $R$
ranges over the regions of the arrangement. In Section
\ref{sec:subdivisions} we show that these images fit together in
an unusual way to give the full image $\Im(g)$; as a consequence,
we obtain two different polyhedral subdivisions of $\Im(g)$, one
of which is a refinement of the other. The finer subdivision is
particularly nice combinatorially; the number of its facets is
$2{m \choose 2}{n \choose 2}r_{m-2,n-2}$, where $r_{m-2,n-2}$ is
the number of regions of the arrangement $\A_{m-2,n-2}$.

This project was suggested in a recent paper of Pachter and
Sturmfels \cite{Pachter}. It arose in their study of the tropical
geometry of statistical models, and provides a new perspective on
Develin, Santos and Sturmfels's study of the image of $g$
\cite{Develin, troprank}.

To explain their motivation, we start with a brief overview of
some recent developments in tropical geometry which led to its
consideration.

\medskip

\noindent 1. The \emph{tropical semiring} $(\R \cup \{\infty\},
\+, \x)$ is the set of real numbers augmented by infinity,
together with the operations of \emph{tropical addition and
multiplication}, which are defined by $x \+ y = \min(x,y)$ and $x
\x y = x+y.$

\emph{Tropical algebraic geometry} is, roughly speaking, the
geometry of the tropical semiring. Let $K = \C\{\{t\}\}$ be the
(algebraically closed) ring of Puiseux series; its elements are
formal power series of the form $f = c_1 t^{a_1} + c_2 t^{a_2} +
\cdots$, where $a_1<a_2<\cdots$ are rational numbers with a common
denominator. There is a natural valuation $\deg: K^* \rightarrow
\Q$, which sends the non-zero Puiseux series $f$ to its degree
$a_1$.

For an ideal $I$ in $K[x_1, \ldots, x_n]$, let $V(I)$ be the
corresponding algebraic variety, intersected with the torus
$(K^*)^n$. It consists of the $n$-tuples $u(t) = (u_1(t), \ldots,
u_n(t))$ of non-zero Puiseux series such that $f(u(t))=0$ for all
$f \in I$. Define $\deg u(t) = (\deg u_1(t), \ldots, \deg
u_n(t))$.

\begin{theorem} \cite{Kapranov, Speyer, Sturmfels1}
For an ideal $I$ in $K[x_1, \ldots, x_n]$, the following subsets
of $\R^n$ coincide.
\begin{enumerate}
\item The topological closure of $\deg V(I)$.

\item The set of $w \in \R^n$ such that the initial ideal
$\iin_w(I)$ contains no monomials.
\end{enumerate}
This set is called the \emph{tropical variety} of $I$, and denoted
$\trop V(I)$.
\end{theorem}

Let $V^+(I)$ be the set of $n$-tuples $u(t)$ in $V(I)$ such that
the coefficient of the leading term of each $u_i(t)$ is a positive
real number.

\begin{theorem} \cite{Lauren}
For an ideal $I$ in $K[x_1, \ldots, x_n]$, the following subsets
of $\R^n$ coincide.
\begin{enumerate}
\item The topological closure of $\deg V^+(I)$.

\item The set of $w \in \R^n$ such that  $\iin_w(I) \cap \R^+[x_1,
\ldots, x_n] = \emptyset$.
\end{enumerate}

This set is called the \emph{positive part of the tropical
variety} of $I$, and denoted $\troppos V(I)$.
\end{theorem}

For all $I$, $\trop V(I)$ is a polyhedral complex, and $\troppos
V(I)$ is a subcomplex of it.

\medskip

\noindent 2. Let $f=(f_1, \ldots, f_n):\R^d \rightarrow \R^n$ be a
polynomial map. Say that $f$ is \emph{positive} if each polynomial
$f_i$ has only positive coefficients. Say that $f$ is
\emph{surjectively positive} if, additionally, every point in the
image of $f$ whose coordinates are positive has a preimage whose
coordinates are positive; \emph{i.e.}, $f(\R_{>0}^d) = \Im f \cap
\R_{>0}^n$.

The \emph{tropicalization of $f$} is the piecewise linear map
$g:\R^d \rightarrow \R^n$ obtained from $f$ by replacing every
$\times$ with a $\x$ and every $+$ with a $\+$. Such a map is
called a \emph{tropical morphism}. \cite{Pachter}

The relationship between the positive part of a tropical variety
and morphisms is outlined in the following theorem.

\begin{theorem} \cite{Pachter, Lauren}
Let $f: \R^d \rightarrow \R^n$ be a polynomial map, and let $g$ be
the tropicalization of $f$. The image $\Im(f)$ is an algebraic
variety in $\R^n$; let $I$ be its corresponding ideal.
\begin{enumerate}
\item $\Im(g) \subset \trop V(I)$.

\item If $f$ is positive, $\Im(g) \subset \troppos V(I)$.

\item If $f$ is surjectively positive, $\Im(g) = \troppos V(I)$.

\end{enumerate}

\end{theorem}

\medskip

\noindent 3. Pachter and Sturmfels \cite{Pachter} used this setup
to study the tropical geometry of statistical models. The
\emph{naive Bayes model with two features} \cite{Garcia} is
parameterized by the
polynomial map $f: \R^{2(m+n)} \rightarrow
\R^{m \times n}$ which maps an $m \times 2$ matrix and a $2 \times
n$ matrix to their product. Its tropicalization, $g$, is the
object of study of this paper.

The image of $f$ consists of the real matrices of rank at most
$2$. This is an algebraic variety, whose corresponding prime ideal
$I$ is generated by the $3 \times 3$ subdeterminants of an $m
\times n$ matrix.

The tropical variety $\trop V(I)$ is the set of $m \times n$ real
matrices having tropical rank at most $2$, and the positive part
$\troppos V(I)$ is the set $B_{m,n}$ of $m \times n$ real matrices
having Barvinok rank at most $2$. For the relevant definitions and
further information, see \cite{Develin,troprank, Pachter}.

This model interests us because the map $f$ is surjectively
positive:

\begin{theorem} \cite{Cohen}
Every positive $m \times n$ matrix of rank $2$ can be written as
the product of a positive $m \times 2$ matrix and a positive $2
\times n$ matrix.
\end{theorem}

It follows that $\Im(g) = \troppos V(I)=B_{m,n}$.

The study of the map $g$ is therefore closely related to Pachter
and Sturmfels's study of the tropical geometry of statistical
models, and to Develin, Santos and Sturmfels's study of the space
$B_{m,n}$.

\section{The morphism and its domains of linearity.} \label{sec:faces}

The piecewise linear map $g: \R^{2(m+n)} \rightarrow \R^{m \times
n}$ which we wish to study is the tropicalization of matrix
multiplication; it is given by
\[
g \left( \begin{bmatrix} & \\ a & A\\ & \end{bmatrix},
\begin{bmatrix} \qquad & b & \qquad \\ \qquad & B & \qquad
\end{bmatrix} \right) = \begin{bmatrix} & \\ a & A\\ &
\end{bmatrix} \x \begin{bmatrix} \qquad & b & \qquad \\ \qquad & B &
\qquad
\end{bmatrix}.
\]
Here $a$ and $A$ denote vectors in $\R^m$ and $b$ and $B$ denote
vectors in $\R^n$. Let $x_i = a_i - A_i$ and $y_j = B_j - b_j$.
The entry $i,j$ of $g$ is given by
\[
g_{ij} = \min(a_i + b_j, A_i + B_j) = \left\{ \begin{matrix} a_i +
b_j, & \mbox{if } x_i \leq y_j \\ A_i + B_j, & \mbox{if } x_i >
y_j. \end{matrix} \right.
\]

The piecewise linear map $g$ is linear on the faces of the
hyperplane arrangement in $\R^{2(m+n)}$:
\[
\A_{m,n}: \qquad x_i = y_j \qquad 1 \leq i \leq m, 1 \leq j \leq
n.
\]
This is essentially the \emph{graphical arrangement} of the
complete bipartite graph $K_{m,n}$; it contains one hyperplane for
each edge of $K_{m,n}$.

Therefore, to understand the fibers of the map $g$, we first study
the faces of $\A_{m,n}$.

Let us briefly outline the close connection between the
arrangement $\A_{m,n}$ and the graph $K_{m,n}$. The matroid
$M(\A_{m,n})$ associated to the arrangement is isomorphic to the
matroid $M(K_{m,n})$ associated to the graph. In other words, a
subset of $k$ hyperplanes of $\A_{m,n}$ intersects in codimension
$k$ if and only if the corresponding subset of $k$ edges of
$K_{m,n}$ does not contain a cycle. Let us call this matroid
simply $M_{m,n}$; it encodes much of the combinatorial information
of the arrangement of the graph. For more details, see \cite[Ch.
2]{Orlik} or \cite[Ch.5]{Oxley}.

The regions of the arrangement $\A_{m,n}$ are in one-to-one
correspondence with the acyclic orientations of the bipartite
graph $K_{m,n}$, as follows: Denote the vertices of $K_{m,n}$ by
$u_1, \ldots, u_m, v_1, \ldots, v_n$, and consider an acyclic
orientation $\o$ of $K_{m,n}$. The corresponding region $R(\o)$
consists of the points $(x,y) \in \R^{m+n}$ such that $x_i < y_j$
if the edge $u_iv_j$ is directed $u_i \rightarrow v_j$ in $\o$,
and $x_i > y_j$ otherwise.

The matroid $M_{m,n}$ appeared independently in the work of Martin
and Reiner \cite{Martin}. They used the finite field method of
\cite{mitesis} to compute a generating function for
$\Chi_{m,n}(q,t)$, the coboundary polynomial of $M_{m,n}$. This
polynomial is a simple transformation of the Tutte polynomial, and
it captures much of the interesting enumerative information of the
matroid.

\begin{proposition} \cite{Martin}
If $\Chi_{m,n}(q,t)$ is the coboundary polynomial of $M_{m,n}$,
\[
1 + q \left( \sum_{(m,n) \in \N^2 - \{(0,0)\}} \Chi_{m,n} (q,t)
\frac{x^m}{m!} \frac{y^n}{n!} \right) = \left( \sum_{(m,n) \in
\N^2} t^{mn} \frac{x^m}{m!} \frac{y^n}{n!} \right)^q.
\]
\end{proposition}

\begin{corollary}\label{orientations} \cite[Ex. 5.6]{EC2}
Let $r_{m,n}$ be the number of acyclic orientations of $K_{m,n}$
(or equivalently, the number of regions of $\A_{m,n}$). Then
\[
\sum_{m,n \geq 1} r_{m,n} \frac{x^m}{m!} \frac{y^n}{n!} =
\frac{e^{x+y}}{e^x+e^y-e^{x+y}}.
\]
\end{corollary}

\begin{proof}
This follows from the formula
$r(m,n)=(-1)^{m+n-1}\Chi_{m,n}(-1,0)$ \cite{Zaslavsky} for the
number of regions of a real arrangement.
\end{proof}

We now proceed to describe and count the domains of linearity of
$g$; that is, the faces of the arrangement $\A_{m,n}$.

\begin{nproposition}\label{faces}
Let $f_{k,m,n}$ be the number of $k$-dimensional faces of the
arrangement $\A_{m,n}$. Then
\[
\sum_{k,m,n \geq \, 0} f_{k,m,n} \, t^k \frac{x^m}{m!} \frac
{y^n}{n!} = \frac{1}{e^{-tx} + e^{-ty} - t(e^x-1)(e^y-1) -1}.
\]
\end{nproposition}

\begin{proof}
To describe a face $F$ of the arrangement we must specify, for
each hyperplane $x_i = y_j$ $(1 \leq i \leq m, \, 1 \leq j \leq
n)$, whether $F$ is in the halfspace $x_i > y_j$, in the halfspace
$x_i<y_j$, or on the hyperplane $x_i=y_j$. Consider the
inequalities or equalities $x_i \bigcirc y_j$ that define a face
$F$. Each one of them puts a restriction on the relative order of
the variables $x_1, \ldots, x_m, y_1, \ldots, y_n$. For a point
$(x,y)$ in $F$, the relative order of $x_i$ and $y_j$ is
determined for all $i$ and $j$. The relative order of $x_{i_1}$
and $x_{i_2}$ is not always determined: an equality $x_{i_1} =
x_{i_2}$ can only be deduced from two defining equalities of the
form $x_{i_1} = y_j$ and $x_{i_2} = y_j$. Similarly, an inequality
$x_{i_1} < x_{i_2}$ can only be deduced as a consequence of two
defining inequalities or equalities of the form $x_{i_1} \bigcirc
y_j$ and $x_{i_2} \bigcirc y_j$.

The faces of $\A_{m,n}$ are described putting these restrictions
together. For example, one face of $\A_{7,8}$ consists of the
points $(x,y) \in \R^{15}$ such that
\[
x_1, x_3 < y_2, y_5, y_7 < x_5 < x_2=y_3=y_6 < y_1 < x_6 < x_4 =
x_7 = y_4 < y_8.
\]

Each face of $\A_{m,n}$ can be described in a similar way, as a
sequence of blocks of variables. The variables $x_{i_1}$ and
$x_{i_2}$ are in the same block if the comparisons $x_{i_1}
\bigcirc y_j$ and $x_{i_2} \bigcirc y_j$ yield the same result for
each $j$. A similar statement holds for $y_{j_1}$ and $y_{j_2}$.
The variables $x_i$ and $y_j$ are in the same block if $x_i = y_j$
in $F$.

Call a block \emph{positive} if it only contains $x_i$s,
\emph{negative} if it only contains $y_j$s, and \emph{mixed}
otherwise. For a point to belong to the face, the variables within
each mixed block are equal. Within an unmixed block, the relative
order of the variables is not determined.

It follows that the faces of $\A_{m,n}$ are in one-to-one
correspondence with the ordered partitions of the set $\{x_1,
\ldots, x_m, y_1, \ldots, y_n\}$ containing no two consecutive
unmixed blocks of the same sign. The dimension of a face is easily
determined from the partition: it is equal to the sum of the sizes
of the unmixed blocks plus the number of mixed blocks. We can now
use the methods of \cite[Chapter 5]{EC2} to compute the desired
generating function.

The generating function for non-empty positive blocks is $X(t,x,y)
= \sum_{n \geq 1} t^n \frac{x^n}{n!} = e^{tx}-1$. The generating
function for non-empty negative blocks is $Y(t,x,y) = e^{ty}-1$.
Therefore, the generating function for partitions of $\{x_1,
\ldots, x_m, y_1, \ldots, y_n\}$ into unmixed blocks of
alternating sign is
\begin{eqnarray*}
Z(t,x,y) &=& (1+X)(1+YX+YXYX+YXYXYX+\cdots)(1+Y) \\ & = &
\frac{1}{e^{-tx} + e^{-ty}-1}.
\end{eqnarray*}

On the other hand, the generating function for mixed blocks is
given by $M(t,x,y) = \sum_{m,n \geq 1} t \frac{x^m}{m!}
\frac{y^n}{n!} = t(e^x-1)(e^y-1)$. The partitions we wish to count
are alternating sequences of partitions of the type counted by $Z$
(which may be empty) and mixed blocks. Therefore

\begin{eqnarray*}
\sum_{k,m,n \geq \, 0} f_{k,m,n} \, t^k \frac{x^m}{m!} \frac
{y^n}{n!} = Z+ZMZ+ZMZMZ+\cdots,
\end{eqnarray*}

which is equal to the given expression.

\end{proof}

Observe that, under the above correspondence, the regions of
$\A_{m,n}$ correspond to the ordered partitions which contain no
mixed blocks. These are counted by $Z(t,x,y)$; setting $t=1$
recovers the generating function for $r_{m,n}$.

It will be convenient to label each region $R$ of the arrangement
with the permutation $\pi(R)$ of the set $[m,\n] = \{1,2,\ldots,m,
\1, \2, \ldots, \n\}$ obtained by reading the blocks from left to
right. The variable $x_i$ corresponds to the letter $i$ (which we
call a \emph{positive letter}), and the variable $y_j$ corresponds
to the letter $\j$ (which we call a \emph{negative letter}).
Within each block, the letters are arranged in increasing order.

For example, the region
\[
y_2, y_4 < x_3 < y_1 < x_1, x_2 < y_3, y_5
\]
of $\A_{3,5}$ is simply denoted by the permutation
$\2\43\112\3\5$.

This labelling is a bijection between the regions of $\A_{m,n}$
and the permutations of $[m,\n]$ such that any two consecutive
letters of the same sign are in increasing order.

\section{The image of $g$ in each region} \label{sec:images}

Consider a region $R$ of $\A_{m,n}$. As we mentioned earlier, the
restriction of the map $g$ to the region $R$ is a linear function.
We now describe the image $g(R)$.

Color each entry of the $m \times n$ matrix $g$ either
\emph{black} or \emph{white}: the entry $g_{ij}$ is black if
$a_i-A_i > B_j-b_j$ in $R$, and white if $a_i-A_i < B_j-b_j$ in
$R$. Permute the rows and columns of $g$ according to the order in
which their labels appear in $\pi$. A path $P$ separates the white
and black entries: it starts at the northwest corner of the
matrix, and takes a step south for each positive letter in $\pi$
and a step east for each negative letter in $\pi$, in the order
prescribed by $\pi$.
\begin{figure}[h]
\centering
\includegraphics[width=4in]{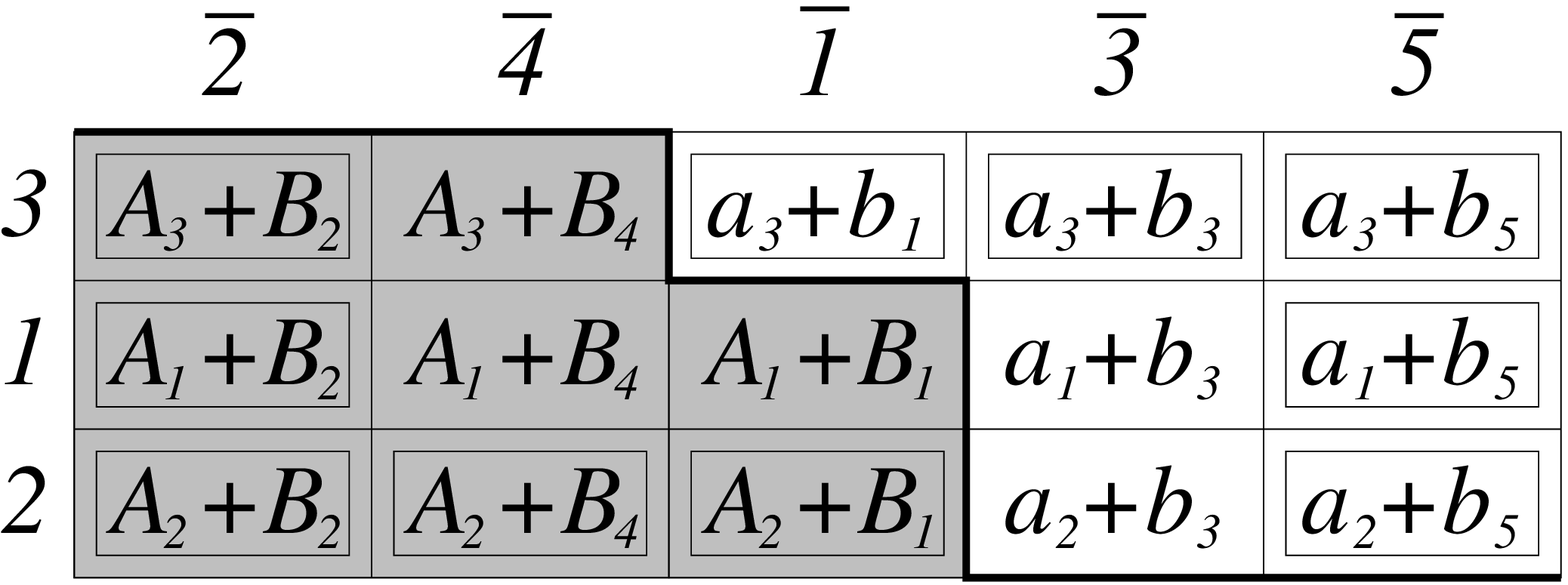}
\caption{The diagram of the region $\2\43\112\3\5$.}\label{matriz}
\end{figure}
Figure \ref{matriz} shows the resulting matrix and path for $\pi =
\2\43\112\3\5$. The entry $g_{ij}$ is $A_i + B_j$ if it is below
$P$ or $a_i+b_j$ if it is above $P$. We call this picture the
\emph{diagram} of $R$ or of $g(R)$.

For $g \in \R^{m \times n}$, write
\[
\Delta_{i_1i_2j_1j_2}(g) = g_{i_1j_1} + g_{i_2j_2} -
g_{i_1j_2}-g_{i_2j_1}.
\]
For simplicity, we will omit $g$ from the notation and simply
write $\Delta_{i_1i_2j_1j_2}$ for this expression. We will call an
equality or inequality of the form $\Delta_{i_1i_2j_1j_2} \bigcirc
\, 0$ a \emph{rectangle relation}. Write $i < \j$ when $i$ appears
before $\j$ in $\pi(R)$.

\begin{nproposition}\label{im(R)1} (Version 1.)
The image $g(R)$ is an open polytope, described by the rectangle
equalities and inequalities that it satisfies. They are the
following:
\begin{eqnarray*}\label{image}
\Delta_{i_1i_2j_1j_2} &=& 0 \textrm{ if } i_1,i_2<\j_1,\j_2, \\
\Delta_{i_1i_2j_1j_2} &=& 0 \textrm{ if } \j_1, \j_2 < i_i, i_2, \\
\Delta_{i_1i_2j_1j_2} &>& 0 \textrm{ if } i_1<\j_1<i_2<\j_2, \\
\Delta_{i_1i_2j_1j_2} &>& 0 \textrm{ if } \j_1 < i_1<\j_2<i_2, \\
\Delta_{i_1i_2j_1j_2} &>& 0 \textrm{ if } i_1 < \j_1<i<\j_2<i_2 \textrm{ for some } i, \\
\Delta_{i_1i_2j_1j_2} &>& 0 \textrm{ if } \j_1 < i_1<\j<i_2<\j_2
\textrm{ for some $\j$,}
\end{eqnarray*}
and no others.
\end{nproposition}

\begin{figure}[h]
\centering
\includegraphics[width=5in]{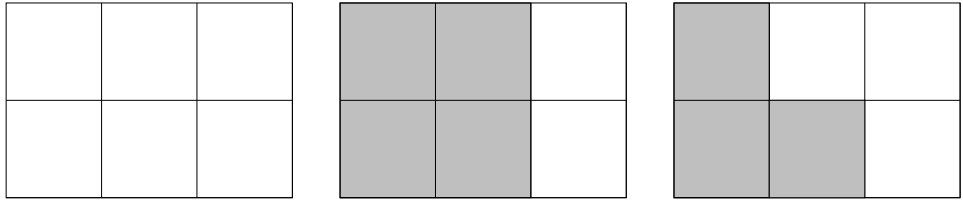}
\caption{The three types of rectangle: monochromatic, sliced, and
jagged.}\label{rectangles}
\end{figure}

In the diagram of $R$, let $\Box_{i_1i_2j_1j_2}$ be the
sub-rectangle of the diagram containing rows $i_1$ through $i_2$
and columns $\j_1$ through $\j_2$. (We implicitly assume that $i_1
< i_2$ and $\j_1 < \j_2$ in $\pi(R)$.) Call this sub-rectangle
\emph{monochromatic} if all its entries have the same color. Call
it \emph{sliced} if it has black and white entries, separated by a
single vertical or horizontal line. Call it \emph{jagged}
otherwise. These definitions are illustrated in Figure \ref{rectangles}.

\addtocounter{nproposition}{-1}

\newpage

\begin{nproposition}\label{im(R)2} (Version 2.)
The image $g(R)$ consists of those $g \in \R^{m \times n}$ such
that $\Delta_{i_1i_2j_1j_2}$ is:
\begin{itemize}
\item equal to zero if $\Box_{i_1i_2j_1j_2}$ is monochromatic,
\item positive if $\Box_{i_1i_2j_1j_2}$ is jagged.
\end{itemize}
\end{nproposition}
\noindent (If $\Box_{i_1i_2j_1j_2}$ is sliced,
$\Delta_{i_1i_2j_1j_2}$ takes positive and negative values in
$R$.)

\begin{proof} It is straightforward to verify that the two versions
of Proposition \ref{im(R)1} are equivalent. Consider a point
$(a,A,b,B) \in \R^{2(m+n)}$ in $R$. The entry $g_{ij}$ of
$g(a,A,b,B)$ is given by $a_i+b_j$ if it is white (\emph{i.e.} if
$i<\j$ in $\pi(R)$), or $A_i+B_j$ if it is black (\emph{i.e.} if
$i>\j$ in $\pi(R)$). These entries are easily seen to satisfy the
given equalities and inequalities.

Conversely, consider a matrix $g \in \R^{m \times n}$ which
satisfies the given equalities and inequalities. Permute its
columns and rows according to the order in which their labels
appear in $\pi$, and draw the path $P$.

Choose appropriate values of the $A_i$s and $B_j$s so that
$g_{ij}=A_i+B_j$ for all the southernmost and westernmost black
entries (the ones surrounded by a box in Figure \ref{matriz}).
This can be done because the system of equations is independent,
and has more unknowns than equations. Make the $A_i$s and $B_j$s
which do not appear in this system of equations very large
positive numbers. Similarly, choose appropriate values of the
$a_i$s and $b_j$s so that $g_{ij}=a_i+b_j$ for all the
northernmost and easternmost white entries, and make the other
$a_i$s and $b_j$s very large positive numbers.

Since $g$ satisfies the rectangle equalities of $g(R)$,
$g_{ij}=A_i+B_j$ for all black entries and $g_{ij}=a_i+b_j$ for
all the white entries.

Take a black entry $g_{ij}$. If the northernmost entry on its
column or the easternmost entry on its row is black, then
$a_i+b_j$ is a very large positive number, larger than $A_i+B_j$.
Otherwise, the rectangle determined by it and the northeasternmost
corner of the matrix has three white entries and one black entry,
and the corresponding rectangle inequality is equivalent to
$A_i+B_j < a_i + b_j$. Similarly, $a_i+b_j < A_i + B_j$ for white
entries $g_{ij}$. Therefore $(a,A,b,B) \in R$ and $g =
g(a,A,b,B)$.

\end{proof}

\begin{corollary}
Let $R$ be a region of $\A_{m,n}$. Then $\span g(R)$ is the
subspace of $\R^{m \times n}$ described by the equalities
\[
g_{i_1j_1} + g_{i_2j_2}= g_{i_1j_2}+g_{i_2j_1}
\]
for those $i_1,i_2,j_1,j_2$ such that $g_{i_1j_1}, g_{i_2j_2},
g_{i_1j_2}$ and $g_{i_2j_1}$ have the same color.
\end{corollary}

For a region $R$ let $\first(R)$ be the length of the first block
of letters of $\pi(R)$ of the same sign (or equivalently, the
number of sources of the corresponding orientation $\o(R)$), and
let $\last(R)$ be the length of the last block of letters of
$\pi(R)$ of the same sign (or equivalently, the number of sinks of
$\o(R)$).

Let $R_1$ be the region where $x_i < y_j$ for all $i$ and $j$, and
let $R_2$ be the region where $x_i > y_j$ for all $i$ and $j$.

\begin{corollary}
Let $R$ be a region of $\A_{m,n}$. If $R \neq R_1, R_2$, then
\[
\dim g(R) = 2m + 2n - 2 - \first(R) - \last(R).
\]
Otherwise, $\dim g(R_1) = \dim g(R_2)= m+n-1$.
\end{corollary}

\begin{proof}
The southernmost and westernmost black entries and the
northernmost and easternmost white entries linearly generate the
remaining ones, and there are no linear relations among them. The
number of them is as claimed.
\end{proof}

In particular, the regions where $g$ has maximum rank are those
corresponding to acyclic orientations of $K_{m,n}$ with a unique
source and a unique sink. This maximum rank is equal to $2m+2n-4$.

%
%
%

\section{Two subdivisions of the image of $g$.}
\label{sec:subdivisions}

The closures of the regions of the hyperplane arrangement
$\A_{m,n}$ give a polyhedral covering of $\R^{2(m+n)}$. Their
images under the map $g$ give a covering of the full image,
$\Im(g)$. We now wish to understand how the images
$\overline{g(R)}$ of the closures of the regions $R$, which we
call the \emph{cells}, fit together.

First notice that if $R$ is a region of $\A_{m,n}$ and $-R$ is its
negative (so the permutation $\pi(-R)$ is equal to the permutation
$\pi(R)$ reversed), then it follows from Proposition \ref{im(R)1}
that $\overline{g(R)} = \overline{g(-R)}$. Therefore, we can
restrict our attention only to the \emph{positive} regions, where
$x_1 < y_1$ (or equivalently, $1 < \1$ in $\pi(R)$).

Also, the following proposition shows that it suffices to study
the maximum-dimensional cells.

\begin{proposition}\cite{Develin}
The image $\Im(g)=B_{m,n}$ is pure.
\end{proposition}

Interestingly, though, our collection of cells \emph{is not}
pure-dimensional. For example, the image of the region
$1\1\2\ldots\n23\ldots m$ is a subspace which is not contained in
any maximum-dimensional cell.

Recall that the cell $\overline{g(R)}$ is maximum-dimensional if
and only if the first and the last block of $\pi(R)$ are
singletons. Call such a maximum-dimensional cell \emph{large} if
the second and the second-to-last blocks of $\pi(R)$ are not
singletons, \emph{small} if the second and the second-to-last
blocks of $\pi(R)$ are both singletons, and \emph{medium} if one
of them is a singleton and the other one is not.

\begin{ntheorem}\label{subdivisions}
The large cells form a polyhedral subdivision of $B_{m,n}$. The
small cells form a finer subdivision of $B_{m,n}$.
\end{ntheorem}

\begin{proof}
The proof is divided into five steps. We start by showing that
each small or medium cell is contained in a large cell; therefore,
the large cells cover $\Im(g)$. Secondly, we show that each large
cell is subdivided into small cells; therefore, the small cells
cover $\Im(g)$ also. The third step is to show that large cells
are pairwise interior-disjoint, and so are small cells. Then we
prove that the covering of $\Im(g)$ with small cells is a
subdivision. Finally, we prove that the covering with large cells
is also a subdivision.

\medskip

\noindent 1. Take any small or medium cell $\overline{g(R_1)}$; we
want to find a large cell containing it. Say $\pi(R_1) =
\j_1i_1\j_2\ldots\j_ri_2\ldots$. Let $R$ be the region with the
label $\pi(R) = i_1\j_1\j_2\ldots\j_ri_2\ldots$: we have adjusted
the beginning of the permutation so that the second block is not a
singleton. This construction is illustrated in Figure
\ref{maximal}. We will show that $\overline{g(R_1)} \subset
\overline{g(R)}$. If $\overline{g(R_1)}$ is medium, then
$\overline{g(R)}$ is large and we are done. If it is small, then
$\overline{g(R)}$ is medium; we can then adjust the end of the
permutation $\pi(R)$ in the same way, to obtain a large cell
containing $\overline{g(R)}$.

\begin{figure}[h]
\centering
\includegraphics[width=5in]{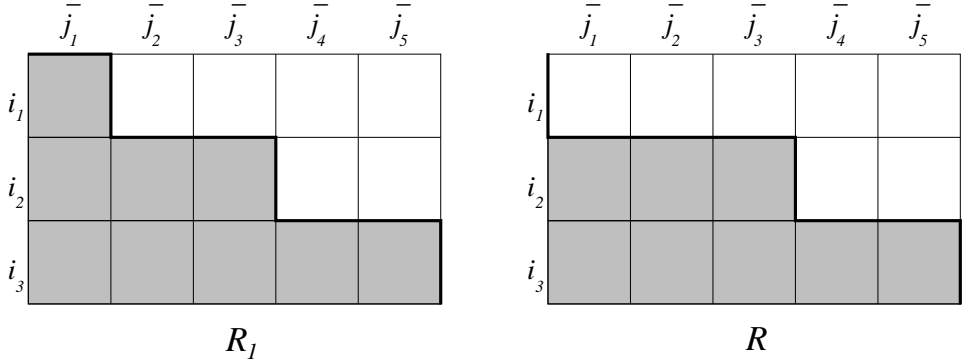}
\caption{A medium cell and its associated large cell.}
\label{maximal}
\end{figure}

We need to show that every rectangle relation satisfied by $g(R)$
is also satisfied by $g(R_1)$. This analysis is most easily
carried out in terms of Version~2 of Proposition \ref{im(R)2}.
Notice that $R$ and $R_1$ have exactly the same monochromatic
sub-rectangles. Also, the only sub-rectangles which are
\emph{sliced} in one diagram and not in the other are
$\Box_{i_1ij_1j_s}$ for $2 \leq s \leq r$ and $i \neq i_1$.

It follows that $g(R_1)$ has the same rectangle relations that
$g(R)$ has, and the additional relations $\Delta_{i_1ij_1j_s}>0$
for $2\leq s \leq r$ and $i \neq i_1$. The rectangle equalities
imply that $\Delta_{i_1ij_1j_s} = \Delta_{i_1i_2j_1j_s}$ for $2
\leq s \leq r$ and $i \neq i_1$. In conclusion, $g(R_1)$ is the
subset of $g(R)$ satisfying the extra relations:
\[
g_{i_1j_1} - g_{i_2j_1} > g_{i_1j_s} - g_{i_2j_s}
\]
for $2 \leq s \leq r$.

We conclude that any maximum-dimensional cell is contained in a
large cell, so the large cells cover $\Im(g)$ by themselves, as
desired.

\medskip

\noindent 2. Now, let us describe which medium and small cells are
contained in a given large cell. Let $R$ be as defined above, and
let $R_s$ be the region with $\pi(R_s) =
\j_si_1\j_1\ldots\widehat{\j_s}\ldots\j_ri_2\ldots$ for $1 \leq s
\leq r$. Imitating the argument above, we conclude that $g(R_s)$
is the subset of $g(R)$ such that $g_{i_1j_s} - g_{i_2j_s}$ is the
unique largest element in $\{g_{i_1j_t} - g_{i_2j_t}: 1 \leq t
\leq r\}$. It follows that $\{\overline{g(R_s)}: 1 \leq s \leq
r\}$ is a subdivision of $\overline{g(R)}$.

In general, the above argument shows that the large cell
$\overline{g(R)}$ is subdivided into $st$ small cells, where $s$
and $t$ are the lengths of the second and second-to-last blocks of
$\pi(R)$. (There are also $s+t$ medium cells in $\overline{g(R)}$.
The first $s$ medium cells form a subdivision of $\overline{g(R)}$
and are subdivided into $t$ small cells each, The other $t$ medium
cells also form a subdivision of $\overline{g(R)}$ and are
subdivided into $s$ small cells each.) In particular, the small
cells also cover $\Im(g)$.

\medskip

\noindent 3. To show that the two coverings of $\Im(g)$ are
polyhedral subdivisions, let us start by showing that the
interiors of any two large cells are disjoint, and the interiors
of any two small cells are also disjoint.

Our strategy is to show that the linear spans of any two large
cells (which are $(2m+2n-4)$-dimensional subspaces) are different;
therefore the intersection of the two cells cannot be
$(2m+2n-4)$-dimensional, and their interiors must be disjoint. We
will use the first Corollary to Proposition \ref{im(R)2}.

Suppose, then, that we know which rectangle equalities an unknown
large cell $\overline{g(R)}$ satisfies. We can recover $R$ as
follows. Define an equivalence relation $\sim$ on $[m] \times [n]$
by declaring that $(i,j) \sim (i',j')$ if $g_{ij}$ and $g_{i'j'}$
are part of the same rectangle equality, and then taking the
transitive closure.

This equivalence relation will have two non-trivial equivalence
classes, and several singletons. (The only exceptions are the
regions with $L$-shaped diagrams, like $i\1\2\ldots
\widehat{\j}\ldots\n 12\ldots \widehat{i} \ldots m \j$. Here there
is only one non-trivial equivalence class, and we can immediately
recover the diagram, and hence the region, from it.)

In one equivalence class, find an entry $(i,j)$ which appears in a
rectangle relation with every other entry in the equivalence
class. Then $g_{ij}$ must be, essentially, the southwest or
northeast corner of the diagram of $R$ which we are trying to
recover. More precisely, $i$ and $\j$ are either in the last
positive and first negative block of $\pi(R)$, or in the last
positive and first negative block, respectively. Because
$\overline{g(R)} = \overline{g(-R)}$, we can assume it is the
former.

Let $\pi(R) = (i_0)\, \j\, (i_1\ldots i_r)\, \sigma \,
(\j_s\ldots\j_1)\, i\, (\j_0)$ be the (unknown) label of $R$. Here
$\sigma$ denotes the segment of $\pi(R)$ which starts at the
second negative letter and ends at the second-to-last positive
letter. Symbols in parenthesis denote letters which may or may not
be in $\pi(R)$.

From Version 1 of Proposition \ref{im(R)1}, $\Delta_{ii'jj'}=0$
holds in $\overline{g(R)}$ if and only if $\j,\j' < i,i'$. For
$i', \j'$ in $\sigma$, this holds if and only if $\j'<i'$. This
allows us to recover which letters are in $\sigma$, and in which
order.

Any positive letters which do not appear in $\sigma$ must appear
to the left of it. Their position in $\pi(R)$ is determined by the
fact that $\overline{g(R)}$ is large: If there is only one such
letter, it must be $i_0$, and $i_1, \ldots, i_r$ do not exist. If
there are several such letters, they must be $i_1, \ldots, i_r$,
and $i_0$ does not exist. Similarly, we recover the positions of
the negative letters which do not appear in $\sigma$. We have
recovered the label of $R$, as claimed.

It follows that the interiors of the large cells are disjoint.
Since each small cell is in a unique large cell, and each large
cell is subdivided into small cells with disjoint interiors, if
follows that the interiors of the small cells are disjoint also.

\medskip

\noindent 4. Let us now show that the covering of $\Im(g)$ by the
small cells is actually a polyhedral subdivision. Consider two
small cells $A_1 = \overline{g(R_1)}$ and $A_2 =
\overline{g(R_2)}$. We wish to show that their intersection is a
face of $A_1$.

A first description of $A_1 \cap A_2$ is given by the rectangle
relations of $A_1$ and $A_2$. We need to find a second
description, which only uses defining equalities and inequalities
of $A_1$ (corresponding to monochromatic and jagged rectangles in
$R_1$), and equalities which define facets of $A_1$ (corresponding
to jagged rectangles in $R_1$).

\begin{figure}[h]
\centering
\includegraphics[width=5in]{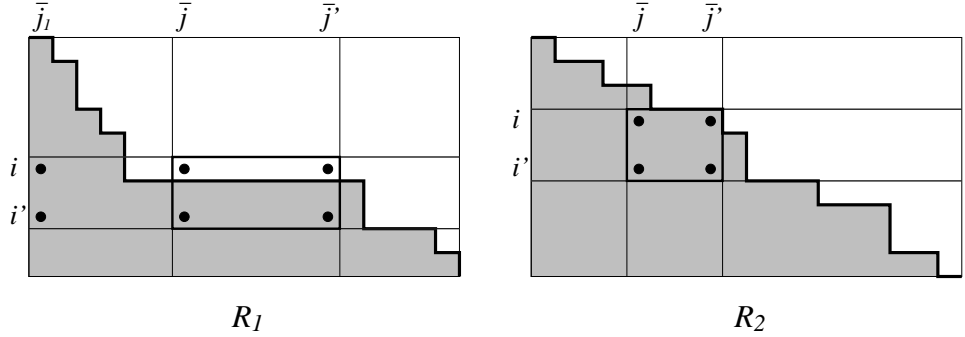}
\caption{The position of $\Box_{ii'jj'}$ in the diagrams of two
small cells.}\label{twosmall}
\end{figure}

Consider a rectangle relation $\Delta_{ii'jj'}=0$ which holds in
$A_1 \cap A_2$. If $\Box_{ii'jj'}$ is monochromatic in $R_1$, the
given relation is satisfied in $A_1$. If it is jagged, then the
relation defines a facet of $A_1$. Therefore, we can assume that
it is sliced. Assume that it is sliced horizontally, so $g_{ij}$
and $g_{ij'}$ are white, and $g_{i'j}$ and $g_{i'j'}$ are black.
The sign of $\Delta_{ii'jj'}$ is undetermined in $A_1$, and
therefore the equality $\Delta_{ii'jj'}=0$ must hold in $A_2$. The
rectangle $\Box_{ii'jj'}$ must then be monochromatic in $R_2$;
assume it is black. Also assume that $i \leq i'$ and $j \leq j'$
in $\pi(R_2)$; our arguments extend immediately to the other
cases. Figure \ref{twosmall} shows the position of $\Box_{ii'jj'}$
in the diagrams of $R_1$ and $R_2$.

Since $A_1$ is small, $\j$ and $\j'$ are not in the first negative
block of $\pi(R_1)$. Let $\j_1$ be the first negative letter in
$\pi(R_1)$. Then $\Box_{ii'j_1 j}$ is jagged in $R_1$, and
$\Delta_{ii'j_1 j} \geq 0$ in $A_1$. If $\j_1 \leq \j$ in
$\pi(R_2)$ then $\Box_{ii'j_1 j}$ is monochromatic in $R_2$ and
$\Delta_{ii'j_1 j} = 0$ in $A_2$; if $\j_1 > \j$ in $\pi(R_2)$
then $\Box_{ii'j j_1}$ is either monochromatic or jagged in $R_2$
and $\Delta_{ii'j j_1} \geq 0$ in $A_2$. In any case,
$\Delta_{ii'j_1j} = 0$ in $A_1 \cap A_2$, and this is a facet
equality of $A_1$. Similarly, $\Delta_{ii'j_1j'} = 0$ in $A_1 \cap
A_2$, and this is a facet equality of $A_1$. Therefore the
equality $\Delta_{ii'jj'}=0$ is a consequence of two facet
equalities of $A_1$.

Exactly the same argument shows that any rectangle inequality
satisfied by $A_1 \cap A_2$ is a consequence of the relations of
$A_1$ and its facet equalities. It follows that the small cells
are actually a polyhedral subdivision of $\Im(g)$.

\medskip

\noindent 5. We now use a similar argument to show that the
covering of $\Im(g)$ into large cells is also a polyhedral
subdivision.

Consider two large cells $A_1 = \overline{g(R_1)}$ and $A_2 =
\overline{g(R_2)}$, and a rectangle relation $\Delta_{ii'jj'}=0$
which holds in $A_1 \cap A_2$. As before, assume that
$\Box_{ii'jj'}$ is sliced horizontally in $R_1$ and monochromatic
black in $R_2$, and that $i \leq i'$ and $\j \leq \j'$ in
$\pi(R_2)$.

The argument that we used for small cells carries over to this
situation, unless $\j$ and $\j'$ are in the first block of
negative letters of $\pi(R_1)$. Since $R_1$ is large, $i$ must
then be the unique first positive letter of $\pi(R_1)$. Figure
\ref{twolarge} shows the position of $\Box_{ii'jj'}$ in the
diagrams of $R_1$ and $R_2$ in this case.
\begin{figure}[h]
\centering
\includegraphics[width=5in]{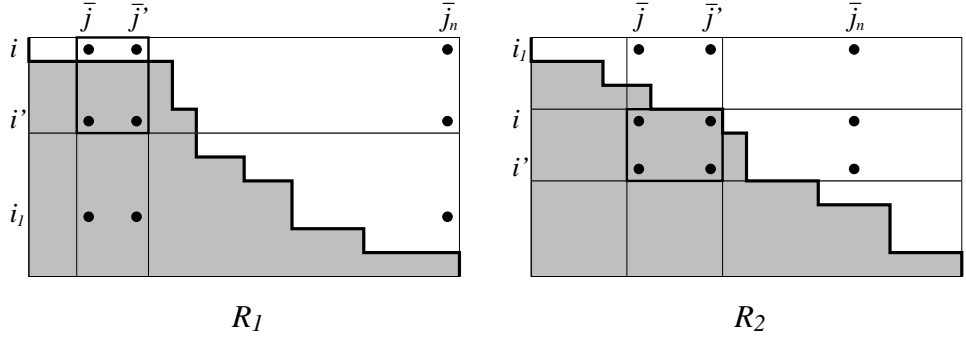}
\caption{The position of $\Box_{ii'jj'}$ in the diagrams of two
large  cells.} \label{twolarge}
\end{figure}

Since $A_1$ is large, $\j$ and $\j'$ are not in the last negative
block of $\pi(R_1)$. Let $\j_n$ be the last negative letter in
$\pi(R_1)$. The rectangles $\Box_{ii'jj_n}$ and $\Box_{ii'j'j_n}$
are jagged in $R_1$. They are either both monochromatic or both
jagged in $R_2$. If they are monochromatic, then
$\Delta_{ii'jj_n}=0$ and $\Delta_{ii'j'j_n}=0$ in $R_2$, and hence
in $R_1 \cap R_2$. These two equalities, which define facets of
$A_1$, will imply the desired equality $\Delta_{ii'jj'}=0$.
Therefore, assume that the two rectangles are jagged in $R_2$,
with $g_{ij_n}$ white.

Now, since $A_2$ is large, $i$ and $i'$ are not in the first
positive block of $\pi(R_2)$. Let $i_1$ be the first positive
letter in $\pi(R_2)$. The rectangle $\Box_{ii_1jj_n}$ is jagged in
$R_1$, so $\Delta_{ii_1jj_n} \geq 0$ in $A_1$. The rectangle
$\Box_{i_1ijj_n}$ is jagged in $R_2$, so $\Delta_{i_1ijj_n} \geq
0$ in $A_2$. Therefore $\Delta_{i_1ijj_n} = 0$ in $A_1 \cap A_2$,
and this is a facet equality of $A_1$. Similarly,
$\Delta_{i_1ij'j_n} = 0$ in $A_1 \cap A_2$, and this is a facet
equality of $A_1$.

It follows that $\Delta_{ii_1jj'} = 0$ in $A_1 \cap A_2$. We also
have that $\Box_{i_1i'jj'}$ is monochromatic in $R_1$ and
$\Delta_{i_1i'jj'} = 0$ holds in $A_1$. The last two equalities
imply that $\Delta_{ii'jj'}=0$ in $A_1 \cap A_2$, and this
equality is implied by two of the facet equalities and one of the
defining equalities of $A_1$. This completes the proof.

\end{proof}

Putting together the considerations in the proofs of Proposition
\ref{im(R)1} and Theorem \ref{subdivisions}, it is now easy to
describe the fiber of a generic point $g$ in the image. Such a
point lies in a unique small cell, in two medium cells, and in a
unique large cell. Therefore, it has preimages in exactly four
regions of $\A_{m,n}$. The diagrams of these four regions are
almost equal; they only differ in the color of the
northwesternmost and southeasternmost entry; assume for simplicity
that they are $g_{11}$ and $g_{mn}$, respectively.

Notice that we can add a constant to the $a_i$s and subtract it
from the $b_j$s, or add a constant to the $A_i$s and subtract it
from the $B_j$s, without affecting $g(a,A,b,B)$. Therefore, we can
focus our attention on the preimages of $g$ with, say, $A_m$ and
$b_n$ fixed.

First consider the preimages in the region where $g_{11}$ and
$g_{mn}$ are black. With $A_m$ and $b_n$ fixed, forcing
$g(a,A,b,B)$ to have the correct border entries determines almost
all of $(a,A,b,B)$. Only $a_m$ and $b_1$ are not determined; $a_m$
can take any value larger than $g_{mn}-b_n$ (the value it would
have if $g_{mn}$ was white), and $b_1$ can take any value larger
than $g_{11}-a_1 = g_{11}-g_{1n}+b_n$ (the value it would have if
$g_{11}$ was white). The preimages form a two-dimensional quadrant
parallel to the $a_mb_1$ plane of $\R^{2(m+n)}$.

The remaining three cells give similar preimages. Putting them all
together, we are left with four two-dimensional quadrants parallel
to the $a_mb_1$, $b_1B_n$, $B_nA_1$ and $A_1a_m$ planes of
$\R^{2(m+n)}$. Their common apex is the point with
$a_m=g_{mn}-b_n, b_1=g_{11}-a_1, B_n=g_{mn}-A_m, A_1=g_{11}-B_1$.

We have two extra degrees of freedom, given by the choices of
$A_m$ and $B_n$. The preimage of a generic point described is a
four-dimensional polyhedral complex; this is consistent with the
dimension drop from the $(2m+2n)$-dimensional range to the
$(2m+2n-4)$-dimensional image.

\medskip

We now have all the tools to enumerate the cells in the two
subdivisions of $B_{m,n}$ that we have constructed.

\begin{nproposition}
Let $\smallcells_{m,n}$ be the number of small cells in $B_{m,n}$.
For $m,n \geq 2$,
\[
\smallcells_{m,n} = 2{m \choose 2}{n \choose 2}r_{m-2,n-2},
\]
where $r_{m-2,n-2}$ is the number of regions of the real
arrangement $\A_{m-2,n-2}$.
\end{nproposition}

\begin{proof}
If we remove the top and bottom rows and leftmost and rightmost
columns from the diagram of a small cell of $B_{m,n}$, we obtain
the diagram of a region of an arrangement combinatorially
isomorphic to $\A_{m-2,n-2}$.

To recover the diagram of the cell from the diagram of the region,
we need to choose the labels of the top and bottom rows (for which
there are $m(m-1)$ options to choose from), and the leftmost and
rightmost columns (for which there are $n(n-1)$ options) which we
deleted. We also need to extend the path that separates the black
and white cells; there is only one way of doing this that gives
the diagram of a small cell. Finally, remember that each cell has
two diagrams that represent it, which differ by a $180^{\circ}$
rotation and a color switch. The desired result follows.
\end{proof}

\begin{nproposition}
Let $\largecells_{m,n}$ be the number of large cells in $B_{m,n}$.
Then
\[
\sum_{m,n \geq 0} \largecells_{m,n}\frac{x^m}{m!}\frac{y^n}{n!} =
\frac12(xX+yY)+ \frac{2XY + X^2(e^x-1) +
Y^2(e^y-1)}{2(e^x+e^y-e^{x+y})},
\]
where $\,X=x(e^y-y-1)$ and $\,Y=y(e^x-x-1)$.

\end{nproposition}

\begin{proof}
Imitate the proof of Proposition \ref{faces}.
\end{proof}

\section{Acknowledgments}

I would like to thank Lior Pachter and Bernd Sturmfels for
suggesting this project, and Mike Develin for helpful
conversations.

\small{

}

\end{document}